\documentclass[11pt]{article}

\usepackage{amsmath,amsfonts,latexsym}
\usepackage{psfrag,epsfig}
\usepackage{young}
\begin{document}

\newtheorem{theorem}{Theorem}[section]
\newtheorem{cor}[theorem]{Corollary}
\newtheorem{lemma}[theorem]{Lemma}
\newtheorem{prop}[theorem]{Proposition}
\newcommand{\sh}{\mbox{ {\rm sh} }}
\newcommand{\bfs}{\mbox{${\bf s}$}}
\newcommand{\Prob}{\mbox{$\mathbb P$}}
\newcommand{\Q}{\mbox{$\mathbb Q$}}
\newcommand{\E}{\mbox{$\mathbb E$}}
\newcommand{\T}{\mbox{$\mathcal T$}}
\newcommand{\bft}{\mbox{$\bf t$}}
\def\enpf{
   {  \parfillskip=0pt\hfil {\hbox{$\Box$}} \par\bigskip  }
   }

\def\R{{\mathbb R}}
\def\C{{\mathbb C}}
\def\Z{{\mathbb Z}}
\def\C{{\mathcal C}}
\def\D{{\mathcal D}}
\def\P{{\mathcal P}}
\def\S{{\mathbb S}}
\def\F{{\cal F}}
\def\i-c{\,{\scriptstyle \bigtriangleup}\,}
\def\s-c{\,{\scriptstyle \bigtriangledown}\,}

\def\l{\lambda}

\begin{center}
\begin{Large}
\begin{bf} 
A path-transformation for random walks
and the Robinson-Schensted correspondence
\end{bf}
\end{Large}

\bigskip

{\sc Neil O'Connell}

\bigskip

{\em \'Ecole Normale Sup\'erieure, Paris}

\end{center}

\bigskip

\begin{abstract}
In [O'Connell and Yor (2002)] a path-transformation $G^{(k)}$ 
was introduced with the property that, for $X$ belonging to a certain 
class of random walks on $\Z_+^k$, the transformed walk $G^{(k)}(X)$
has the same law as that of the original walk conditioned never
to exit the Weyl chamber $\{x:\ x_1\le\cdots\le x_k\}$.  In this paper,
we show that $G^{(k)}$ is closely related to the Robinson-Schensted
algorithm, and use this connection to give a new proof of the above
representation theorem.  The new proof is valid for a larger class
of random walks and yields additional information about the joint
law of $X$ and $G^{(k)}(X)$.  The corresponding results for the
Brownian model are recovered by Donsker's theorem.  These are
connected with Hermitian Brownian motion and the Gaussian Unitary
Ensemble of random matrix theory.  The connection we make between
the path-transformation $G^{(k)}$ and the Robinson-Schensted algorithm also
provides a new formula and interpretation for the latter.
This can be used to study properties of the Robinson-Schensted algorithm and,
moreover, extends easily to a continuous setting.
\end{abstract}

\vfill



\newpage

\section{Introduction and summary}

For $k\ge 2$, denote the set of probability distributions on $\{1,\ldots,k\}$
by $\P_k$.  Let $(\xi_m,\ m\ge 1)$ be a sequence of independent random
variables with common distribution $p\in\P_k$ and, for $1\le i\le k$,
$n\ge 0$, set
\begin{equation}
X_i(n)=|\{1\le m\le n:\ \xi_m=i\}| .
\end{equation}

If $p_1<\cdots <p_k$, there is a positive probability 
that the random walk $X=(X_1,\ldots,X_k)$ never exits the Weyl chamber
\begin{equation}\label{weyl}
W=\{x\in\R^k:\ x_1\le\cdots\le x_k\};
\end{equation}
this is easily verified using, for example, the concentration inequality
$$P(|X(n)-np|>\varepsilon )\le K e^{-c(\varepsilon)n},$$
where $c(\varepsilon)$ and $K$ are finite positive constants.

In~\cite{oy2}, a certain path-transformation $G^{(k)}$ was introduced
with the property that:
\begin{theorem}\label{oy2} The law of the transformed walk $G^{(k)}(X)$ 
is the same, assuming $p_1<\cdots <p_k$, as that of the original
walk $X$ conditioned never to exit $W$.  
\end{theorem}
We will recall the definition of $G^{(k)}$ in Section 2 below.

This was motivated by a
desire to find a multi-dimensional generalisation of Pitman's 
representation for the three-dimensional Bessel process~\cite{pitman},
and to understand some striking connections which were
recently discovered by Baik, Deift and Johansson~\cite{bdj}, 
Baryshnikov~\cite{bar} and Gravner, Tracy and Widom~\cite{gtw},
between oriented percolation and random matrices.  For more
background on this, see~\cite{survey}.

The proof of Theorem~\ref{oy2} given in~\cite{oy2} uses certain
symmetry and reversibility properties of M/M/1 queues in series;
consequently, the transformation $G^{(k)}$ has a `queueing-theoretic'
interpretation.

In this paper we will show that the path-transformation $G^{(k)}$
is closely related to the Robinson-Schensted correspondence.
More precisely, if $\lambda(n)=(\lambda_1(n)\ge\cdots\ge\lambda_k(n))$ 
denotes the shape of the Young tableaux obtained, when one applies
the Robinson-Schensted algorithm, with column-insertion, to the random word
$\xi_1\cdots\xi_n$, then (for any realisation of $X$)
$$(G^{(k)}(X))(n)=(\lambda_k(n),\ldots,\lambda_1(n)).$$

Immediately, this yields a new representation and formula for the 
Robinson-Schensted algorithm, and this formula has a queueing interpretation.
We will use this representation to recover known, and perhaps 
not-so-well-known, properties of the Robinson-Schensted algorithm.

Given this connection, Theorem~\ref{oy2} can now be interpreted 
as a statement about the evolution of the shape $\lambda(n)$ of a 
certain randomly growing Young tableau.  We give a direct proof of 
this result using properties of the Robinson-Schensted correspondence.  This also yields 
more information about the joint law of $X$ and $G^{(k)}(X)$, and 
dispenses with the condition $p_1<\cdots <p_k$.  

As in~\cite{oy2}, the corresponding results for the Brownian motion
model can be recovered by Donsker's theorem.  The path-transformation
$G^{(k)}$ extends naturally to a continuous setting and, given the
connection with the Robinson-Schensted algorithm, the continuous version can now
be regarded as a natural extension of the Robinson-Schensted algorithm to a continuous
setting.  As discussed in~\cite{oy2}, the results for Brownian motion
have an interpretation in random matrix theory.  In particular, 
Theorem~\ref{oy2} yields a representation for the 
eigenvalue process associated with Hermitian Brownian motion as a certain 
path-transformation 
(the continuous analogue of $G^{(k)}$) applied to a standard Brownian motion.
The new results presented in this paper also yield new results in this
context.
This random matrix connection comes from the well-known fact that the 
eigenvalue process 
associated with Hermitian Brownian motion can be interpreted as Brownian
motion conditioned never to exit the Weyl chamber $W$.  We remark that
a similar representation for the eigenvalues of Hermitian Brownian motion
was independently obtained by Bougerol and Jeulin~\cite{bj},
in a more general context, by completely different methods.

The outline of the paper is as follows.
In the next section we recall the definition of $G^{(k)}$
and record some of its properties.  In section 3, we make
the connection with the Robinson-Schensted algorithm, and briefly consider
some immediate implications of this connection.  A worked example
is presented in section 4.  In section 5,
we record some properties of the conditioned walk of Theorem~\ref{oy2}
and extend its definition beyond the case $p_1<\cdots <p_k$.
In section 6, we prove a generalisation of Theorem~\ref{oy2},
in the context of Young tableaux, using properties of the RS
correspondence.  In section 7, we define a continuous version
of the path-transformation and present the `Poissonized' analogues
of the results of the previous section.  In section 8, we present
the corresponding results for the Brownian model, and briefly
discuss the connection with random matrices.  An application
in queueing theory is presented in section 9, and we conclude
the paper with some remarks in section 10.  

{\em Some notation:}
Let $b=\{e_1,\ldots,e_k\}$ denote the standard basis elements
in $\R^k$.
For $x,y\in\R_+^k$ we will write $x^y=x_1^{y_1}\cdots x_k^{y_k}$,
$xy=(x_1y_1,\ldots,x_ky_k)$, $|x|=\sum_ix_i$ and define 
$x^*\in\R_+^k$ by $x^*_i=x_{k-i+1}$. Denote the origin in $\R^k$
by $o$.

{\em Acknowledgements:}
Thanks to Francois Baccelli, Phillipe Biane, Phillipe Bougerol and 
Marc Yor for many helpful and illuminating discussions on these topics.
This research was partly carried out during a visit, funded by the {\em CNRS},
to the Laboratoire de Probabilit\'es, Paris 6, and partly at the {\em ENS},
thanks to financial support of {\em INRIA}.

\section{The path-transformation}

The support of the random walk $X$, which we denote by $\Pi_k$,
consists of paths $x:\Z_+\to\Z_+^k$ with $x(0)=0$ and, for each $n>0$, 
$x(n)-x(n-1)\in b$.  
Let $\Pi_k^W$ denote the subset of those paths taking values in $W$.
It is convenient to introduce another set $\Lambda_k$ of paths
$x:\Z_+\to\Z_+^k$ with $x(0)=0$ and $x(n)-x(n-1)\in\{0, e_1,\ldots,e_k\}$,
for each $n>0$.

For $x,y\in\Lambda_1$, define $x\i-c y\in \Lambda_1$ and 
$x\s-c y\in \Lambda_1$ by
\begin{equation}\label{ic}
(x\i-c y)(n)=\min_{0\le m\le n}[x(m)+y(n)-y(m)] ,
\end{equation}
and
\begin{equation}\label{sc}
(x\s-c y)(n)=\max_{0\le m\le n}[x(m)+y(n)-y(m)].
\end{equation}
The operations $\i-c$ and $\s-c$ are not associative in general.
Unless otherwise deleniated by parentheses, the default order of 
operations is from left to right; for example, when we write $x\i-c y\i-c z$, 
we mean $(x\i-c y)\i-c z$.

The mappings $G^{(k)}: \Lambda_k\to\Lambda_k$ are defined as follows.
Set
\begin{equation}\label{g2}
G^{(2)}(x,y)=(x\i-c y, y\s-c x)
\end{equation} 
and, for $k>2$, 
\begin{align}\label{gn}
G^{(k)}(x_1, &  \ldots,x_k)=(x_1\i-c x_2\i-c\cdots\i-c x_k,\nonumber\\
& G^{(k-1)}(x_2\s-c x_1,x_3\s-c (x_1\i-c x_2),\ldots,x_k\s-c (x_1\i-c
\cdots\i-c x_{k-1}))) .
\end{align}

Note that $G^{(k)}:\Pi_k\to\Pi_k^W$.

We will now give an alternative definition of $G^{(k)}$ which will
be useful for making the connection with the Robinson-Schensted correspondence.

Occasionally, we will suppress the dependence of functions on $x$,
when the context is clear:  for example, we may write
$G^{(k)}$ instead of $G^{(k)}(x)$, and so on.

For $k\ge 2$, define maps $D^{(k)}:\Lambda_k\to \Lambda_k$ and 
$T^{(k)}:\Lambda_k\to \Lambda_{k-1}$, by
\begin{equation}
D^{(k)}(x)=(x_1,x_1\i-c x_2, \ldots, x_1\i-c\cdots\i-c x_k)
\end{equation}
and
\begin{equation}
T^{(k)}(x)=(x_2\s-c x_1,x_3\s-c (x_1\i-c x_2),\ldots,x_k\s-c 
(x_1\i-c\cdots\i-c x_{k-1})) .
\end{equation}
For notational convenience, let $D^{(1)}$
be the identity transformation.  

Note that the above definition is recursive:
for $i\ge 2$,
\begin{equation}
D_i^{(k)}=D^{(k)}_{i-1}\i-c x_i ,
\end{equation}
and 
\begin{equation}
T^{(k)}_{i-1}=x_i\s-c D^{(k)}_{i-1}.
\end{equation}
Alternatively, we can write
\begin{equation}
(D_i^{(k)},T^{(k)}_{i-1})=G^{(2)}(D^{(k)}_{i-1}, x_i).
\end{equation}

For each $x\in\Lambda_k$, consider the triangular array of sequences 
$d^{(i)}\in\Lambda_{k-i+1}$, $1\le i\le k$, defined as follows.
Set $$d^{(1)}=D^{(k)}(x),\ \ \  t^{(1)}=T^{(k)}(x),$$
$$d^{(2)}=D^{(k-1)}(t^{(1)}),\ \ \ t^{(2)}=T^{(k-1)}(t^{(1)}),$$ and
so on; for $i\le k$, $$d^{(i)}=D^{(k-i+1)}(t^{(i-1)}),$$ and
for $i\le k-1$, $$t^{(i)}=T^{(k-i+1)}(t^{(i-1)}).$$ 
Recalling the definition of $G^{(k)}$ given earlier,
we see that
\begin{equation}
G^{(k)}=(d^{(1)}_k,\ldots,d^{(k)}_1).
\end{equation}
Note also that, for each $i\le k$,
\begin{equation}
G^{(i)}(x_1,\ldots,x_i)=(d^{(1)}_i,\ldots,d^{(i)}_1).
\end{equation}

We will conclude this section by recording some useful
properties and interpretations of the operations $\i-c$ and $\s-c$, 
and of the path-transformation $G^{(k)}$, for later reference.
We defer the proofs: these will be given in the appendix.

The following notation for increments of paths will be useful:
for $x\in\Lambda_k$ and $l\ge n$, set $x(n,l)=x(l)-x(n)$.

The operations $\i-c$ and $\s-c$ have a queueing-theoretic
interpretation, which we will make strong use of when we
make the connection with the Robinson-Schensted correspondence in the next
section.  For more general discussions on `min-plus algebra'
and queueing networks, see~\cite{bb}.

Suppose $(x,y)\in\Pi_2$, and consider a simple
queue which evolves as follows.  At each time $n$, either
$x(n)-x(n-1)=1$ and $y(n)-y(n-1)=0$, in which case a new customer
arrives at the queue, or $x(n)-x(n-1)=0$ and $y(n)-y(n-1)=1$,
in which case, if the queue is not empty, a customer departs
(otherwise nothing happens).  The number of customers remaining
in the queue at time $n$, which we denote by $q(n)$, satisfies
the Lindley recursion
\begin{equation}\label{lindley}
q(n)=\max\{ q(n-1)+\epsilon(n), 0\} ,
\end{equation}
where $\epsilon(n)=x(n)-x(n-1)-y(n)+y(n-1)$.
Iterating (\ref{lindley}), we obtain
\begin{equation}\label{q}
q(n)=\max_{0\le m\le n} [ x(m,n)-y(m,n) ] .
\end{equation}
Thus, the number of customers $d(n)$ to depart up to and including
time $n$ is given by
\begin{equation}
d(n)=x(n)-q(n)= (x\i-c y)(n).
\end{equation}
We also have
\begin{equation}
t(n):=x(n)+u(n)=(y\s-c x)(n),
\end{equation}
where $$u(n)=y(n)-d(n)$$ is
the number of times $m\le n$ that $y(m)-y(m-1)=1$ and $q(m-1)=0$;
in the language of queueing theory, $u(n)$ is the number of
`unused services' up to and including time $n$. (For this queue
we refer to the points of increase of $y$ as `services'.)

\begin{lemma}\label{prop1}
For $(x,y)\in\Lambda_2$,
\begin{equation}
x\i-c y+y\s-c x =x+y ,
\end{equation}
and
\begin{eqnarray*}
x(n)-(x\i-c y)(n) &=& \max_{0\le m\le n}[x(m,n)-y(m,n)] \\
&=& \max_{l\ge n} [(x\i-c y)(n,l)-(y\s-c x)(n,l) ].
\end{eqnarray*}
In particular, writing $G^{(2)}\equiv G^{(2)}(x,y)$, we have:
\begin{equation}
(x(n),y(n))=G^{(2)}(n) + F^{(2)}\left( G^{(2)}(n,l),\ l\ge n\right) ,
\end{equation}
where $F^{(2)}:\D\to\Z^2$ is defined, on
$$\D = \{ z\in (\Z^2)^{\Z_+}:\ M(z)=\max_{n\ge 0}[z_1(n)-z_2(n)]<\infty\} ,$$
by $F^{(2)}(z) = (M(z),-M(z))$.
\end{lemma}

In the queueing context described above, the above lemma states
that $x+y=d+t$ and
\begin{equation}\label{fut}
q(n)=\max_{m\ge n} [ d(n,m)-t(n,m) ].
\end{equation}
The first identity is readily verified.
The formula for $q(n)$ in terms of the future increments of $d$
and $t$ follows from the time-reversal symmetry in the dynamics
of the system: this formula is the dual of (\ref{q}).  When time is
reversed, the roles played by $(x,y)$ and $(d,t)$ are interchanged.
This symmetry is at the heart of the proof of Theorem~\ref{oy2}
given in~\cite{oy2}, where it is considered in an equilibrium context.

Note that, if we set $z=y-x$ and $s(n)=\max_{0\le m\le n}z(m)$,
then $$y\s-c x-x\i-c y=2s-z$$ and (\ref{fut}) is equivalent to
the well-known identity $$s(n)=\min_{l\ge n}[2s(l)-z(l)].$$
This is familiar in the context of Pitman's representation
for the three-dimensional Bessel process.  Observe that the
statement of Theorem~\ref{oy2} in the case $k=2$ is equivalent
to the following discrete version of Pitman's theorem:  if
$\{Z(n),n\ge 0\}$ is a simple random walk on $\Z$ with positive
drift, started at $0$, and we set $S(n)=\max_{0\le m\le n}Z(n)$,
then $2S-Z$ has the same law as that of $Z$ conditioned to stay
non-negative.  The usual statement of Pitman's theorem can be
recovered from Theorem~\ref{oy2bm} below.

Lemma~\ref{prop1} has the following generalisation:
\begin{lemma} \label{prop2} 
For $x\in\Lambda_k$, writing $G^{(k)}\equiv G^{(k)}(x)$,
\begin{equation}
|G^{(k)}|=|x| ,
\end{equation}
and
\begin{equation}
x(n)=G^{(k)}(n)+F^{(k)}\left( G^{(k)}(n,l),\ l\ge n\right) ,
\end{equation}
where the function $F^{(k)}$ will be defined in the proof.
\end{lemma}

As we remarked earlier, the operations $\i-c$ and $\s-c$ are not
associative.  The following identities are useful for manipulating
complex combinations of these operations.
\begin{lemma}\label{prop3}
For $(a,b,c)\in\Lambda_3$,
\begin{equation}\label{uf1}
a\s-c (c\i-c b)\s-c (b\s-c c) = a\s-c b \s-c c,
\end{equation}
and 
\begin{equation}\label{ouf1}
a\i-c (c\s-c b)\i-c (b\i-c c) = a\i-c b\i-c c .
\end{equation}
\end{lemma}

For example, (\ref{uf1}) immediately yields:
\begin{lemma}\label{prop4}
For $x\in\Lambda_k$, $G^{(k)}_k(x)=x_k\s-c\cdots\s-c x_1$.
\end{lemma}

\section{Connection with the Robinson-Schensted algorithm}\label{gkrs}

We refer the reader to the books of Fulton~\cite{fulton} and
Stanley~\cite{stanley} for detailed discussions on the Robinson-Schensted 
algorithm and its properties.  The standard Robinson-Schensted algorithm
takes a word $w=a_1\cdots a_n \in \{1,2,\ldots,k\}^n$ and proceeds, 
by `row-inserting' the numbers $a_1$, then $a_2$, and so on, to 
construct a semistandard tableau $P(w)$ associated with $w$,
of size $n$ with entries from the set $\{1,2,\ldots,k\}$.  
If one also maintains a `recording tableau' $Q(w)$, which is
standard tableau of size $n$, the mapping from words 
$\{1,2,\ldots,k\}^n$ to pairs of semistandard and standard
tableaux of size $n$, the semistandard tableau having entries
from $\{1,2,\ldots,k\}$ and both having the same shape, 
is a bijection:  this is the Robinson-Schensted correspondence.  
One can also do all of the above using
`column-insertion' instead of row-insertion to construct
the semistandard tableau, but still maintaining a recording
tableau, and the resulting map is also a bijection.  Column
and row insertion are not the same thing, but they are related
in the following way:  the semistandard tableau obtained by
applying the Robinson-Schensted algorithm, with column-insertion, to the
word $a_1\cdots a_n$ is the same as the one obtained by
applying the Robinson-Schensted algorithm, with row-insertion, to the
reversed word $a_n\cdots a_1$.  The standard tableaux obtained
in each case are also related, but we do not need this
and refer the reader to~\cite{fulton} for details.

Fix $x\in\Pi_k$, the let $d^{(i)}\in\Lambda_{k-i+1}$, $1\le i\le k$,
be the corresponding triangular array of sequences defined in the
previous section.

For each $n$, construct a semistandard Young 
tableau as follows.  In the first row, put 
$$d^{(1)}_1(n)\mbox{  1's,  }d^{(2)}_1(n)-d^{(1)}_1(n)
\mbox{  2's,  }\ldots d^{(k)}_1(n) -d^{(k-1)}_1(n)
\mbox{  $k$'s;  }$$
in the second row, put 
$$d^{(1)}_2(n)\mbox{  2's,  }d^{(2)}_2(n)-d^{(1)}_2(n)
\mbox{  3's,  }\ldots d^{(k-1)}_2(n)-d^{(k-2)}_2(n)
\mbox{  $k$'s;  }$$
and so on.  In the final row, there are just $d^{(1)}_k(n)$
$k$'s.  Denote this tableau by $\tau(n)$.
For example, if $k=3$ and
\begin{equation} 
\begin{array}{ccc}
d^{(1)}_1(7) & d^{(1)}_2(7) & d^{(1)}_3(7) \\
& d^{(2)}_1(7) & d^{(2)}_2(7) \\
& & d^{(3)}_1(7) 
\end{array}
=
\begin{array}{ccc}
2 & 2&1\\
&3&2\\
&&4
\end{array}
\end{equation}
then the corresponding semistandard tableau $\tau(7)$ is
\begin{equation}
\begin{Young}
1&1&2&3\cr
2&2\cr
3\cr
\end{Young}
\end{equation}

Let $a_m$ be the sequence defined by $a_m=i$ whenever
$$x(m)-x(m-1)=e_i.$$

\begin{theorem}\label{tgkrs}
The semistandard tableau $\tau(n)$ is precisely the one which is
obtained when one applies the Robinson-Schensted algorithm, with column insertion, 
to the word $a_1\cdots a_n$.  In particular, 
if $l(n)$ denotes the shape of $\tau(n)$, then
$l(n)^*=(G^{(k)}(x))(n)$.
\end{theorem}

{\em Proof:}  It will suffice to describe how the mapping $G^{(k)}$ acts
on a typical element of $\Pi_k$, from an algorithmic point of view.

For each $k\ge 2$, the maps $D^{(k)}:\Lambda_k\to \Lambda_k$ and 
$T^{(k)}:\Lambda_k\to \Lambda_{k-1}$ can be defined as follows.  
Fix $x\in\Pi_k$ and set $d=D^{(k)}(x)$, $t=T^{(k)}(x)$.  Set $d(0)=t(0)=0$,
and define the sequences $d(n)$ and $t(n)$ inductively on $n$.
Suppose $x(n)-x(n-1)=e_i$; that is, $a_n=i$.  
We need to treat the cases $i=1$ and $i=k$ separately.  

Suppose $i=1$.  Then we set $d(n)=d(n-1)+e_1$
and $t(n)=t(n-1)+e_1$.  

If $i=k$, and $d_k(n-1)<d_{k-1}(n-1)$,
we set $d(n)=d(n-1)+e_k$ and $t(n)=t(n-1)$.  

If $i=k$, and
$d_k(n-1)=d_{k-1}(n-1)$, we set $d(n)=d(n-1)$ and $t(n)=t(n-1)+e_{k-1}$.

Now suppose $1<i<k$.  If $d_i(n-1)<d_{i-1}(n-1)$, set $d(n)=d(n-1)+e_i$
and $t(n)=t(n-1)+e_i$; if $d_i(n-1)=d_{i-1}(n-1)$, set $d(n)=d(n-1)$
and $t(n)=t(n-1)+e_{i-1}$.  Recall that $D^{(1)}$
is the identity transformation.

In queueing language, we have just constucted a series of $k$ queues
in tandem. Initially there are infinitely many customers in the
first queue and the other queues are all empty.  At each time $n$,
if $x_i(n+1)-x_i(n)=1$ (or, equivalently, $a_n=i$)
there is a `service event' at the $i^{th}$ queue;
if this queue is not empty a customer departs from it and,
if $i<k$, joins the $(i+1)^{th}$ queue.  The number
of departures from the $i^{th}$ queue up to and including time $n$
is given by $d_i(n)$ and $t_i(n)=d_i(n)+u_i(n)$, where $u_i(n)$
is the number of `unused services' at the $(i+1)^{th}$ queue up to
and including time $n$.  Recalling the queueing-theoretic interpretations 
of $\i-c$ and $\s-c$, we see that $d_i=x_1\i-c\cdots\i-c x_i$ 
and $t_i=x_i\s-c d_{i-1}$.

Now fix $x\in\Pi_k$, and recall the definition of the
triangular array of sequences $d^{(i)}\in\Lambda_{k-i+1}$, $1\le i\le k$.
Set $$d^{(1)}=D^{(k)}(x),\ \ \  t^{(1)}=T^{(k)}(x),$$
$$d^{(2)}=D^{(k-1)}(t^{(1)}),\ \ \ t^{(2)}=T^{(k-1)}(t^{(1)}),$$ and
so on; for $i\le k$, $$d^{(i)}=D^{(k-i+1)}(t^{(i-1)}),$$ and
for $i\le k-1$, $$t^{(i)}=T^{(k-i+1)}(t^{(i-1)}).$$ 
Here we have constructed a `series of queues in series', the entire
system `driven' by $x$.  This is represented in Figure 1 for the
case $k=3$.  

\begin{figure}
\setlength{\unitlength}{1cm}
\begin{picture}(12,5)
\begin{small}
\psfrag{x1}{$x_1$}
\psfrag{x2}{$x_2$}
\psfrag{x3}{$x_3$}
\psfrag{xk}{$x_k$}
\psfrag{q11}{$q^{(1)}_1$}
\psfrag{d11}{$d^{(1)}_1$}
\psfrag{q12}{$q^{(1)}_2$}
\psfrag{d12}{$d^{(1)}_2$}
\psfrag{d13}{$d^{(1)}_3$}
\psfrag{d1k}{$d^{(1)}_k$}
\psfrag{q21}{$q^{(2)}_1$}
\psfrag{d21}{$d^{(2)}_1$}
\psfrag{d22}{$d^{(2)}_2$}
\psfrag{d1k-1}{$d^{(1)}_{k-1}$}
\psfrag{d2k-1}{$d^{(2)}_{k-1}$}
\psfrag{d2k-2}{$d^{(2)}_{k-2}$}
\psfrag{d31}{$d^{(3)}_1$}
\psfrag{dk-12}{$d^{(k-1)}_2$}
\psfrag{dk1}{$d^{(k)}_1$}
\psfrag{t11}{$t^{(1)}_1$}
\psfrag{t12}{$t^{(1)}_2$}
\psfrag{t21}{$t^{(2)}_1$}
\psfrag{t1k-1}{$t^{(1)}_{k-1}$}
\psfrag{tk-11}{$t^{(k-1)}_1$}
\psfrag{tk-22}{$t^{(k-2)}_2$}
\put(0,0){\epsfig{file=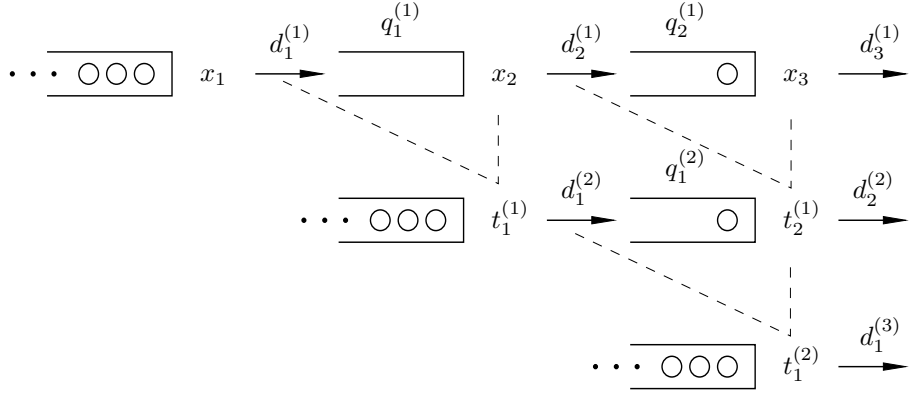,height=5cm,width=12cm}}
\end{small}
\end{picture}
\caption{The series of queues in series ($k=3$)}
\end{figure}

The first series of queues is just the one described above:
there are $k$ queues in series, and initially queues $2$
thru $k$ are empty and the first queue has infinitely many
customers;  whenever $x_i$ increases by one there is a service
at the $i^{th}$ queue and one customer is permitted to depart
(and proceed to the next queue if $i<k$).  The number of
departures from the $i^{th}$ queue up to time $n$ is given by
$d^{(1)}_i(n)$.  

The second series of queues has $t^{(1)}$ `moving' the customers
in place of $x$.  This time there are $k-1$ queues.  
Initially, queues $2$ thru $k-1$ are empty and the first queue 
has infinitely many customers.  There is
a service event at the $i^{th}$ queue whenever $t^{(1)}_i$
increases by one---that is, whenever, in the first series,
there is a departure from the $i^{th}$ queue or an unused
service at the $(i+1)^{th}$ queue (these events will never
occur simultaneously).  

The second series generates a new sequence of `t's', which
we denote by $t^{(2)}$, and this is used to drive the third
series, which consists of $k-2$ queues, and so on.

It is useful to define 
\begin{equation}
q^{(j)}_i = d^{(j)}_i - d^{(j)}_{i+1} ;
\end{equation}
$q^{(j)}_i(n)$ is just the number of customers in the
$(i+1)^{th}$ queue of the $j^{th}$ series at time $n$.
For example, in the network shown in Figure 1, $q^{(1)}_1=0$
and $q^{(1)}_2=q^{(2)}_1=1$.

Now consider the evolution of the 
corresponding semistandard tableaux $\tau(n)$, $n\ge 1$.
See Figure 2.
\begin{figure}
\setlength{\unitlength}{1cm}
\begin{picture}(10,12)
\begin{small}
\psfrag{1n}{1}
\psfrag{2n}{2}
\psfrag{3n}{3}
\psfrag{q11}{$q^{(1)}_1$}
\psfrag{d11}{$d^{(1)}_1$}
\psfrag{q12}{$q^{(1)}_2$}
\psfrag{d12}{$d^{(1)}_2$}
\psfrag{d13}{$d^{(1)}_3$}
\psfrag{q21}{$q^{(2)}_1$}
\psfrag{d21}{$d^{(2)}_1$}
\psfrag{d22}{$d^{(2)}_2$}
\psfrag{d31}{$d^{(3)}_1$}
\put(3.5,2){\epsfig{file=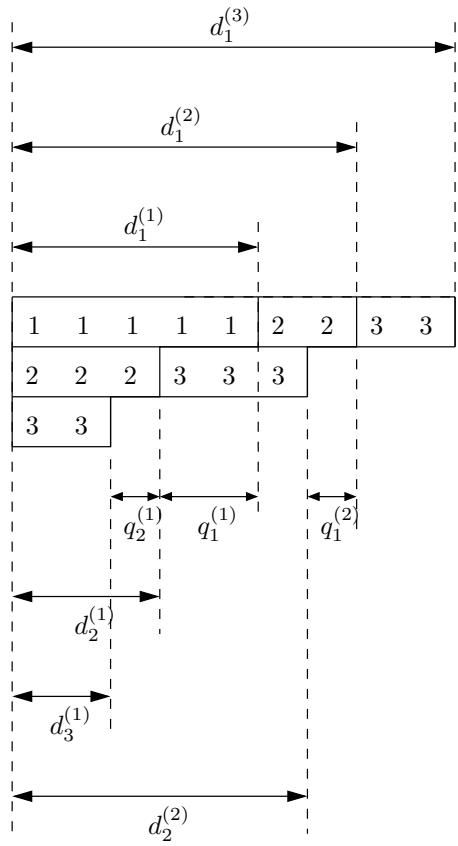,height=11cm,width=6cm}}
\end{small}
\end{picture}
\caption{The tableau $\tau(17)$}
\end{figure}
Another look at the algorithm described above should convince
the reader that $\tau(n)$ is precisely the semistandard
tableau obtained when one applies the Robinson-Schensted algorithm,  with
{\em column} insertion, to the word $a_1\cdots a_n$.

To see this, look at the tableau $\tau(17)$ represented in Figure 2.
Recall that $a_m=i$ if $x_i$ increases by one at time $m$,
in which case there is a service at the $i^{th}$ queue of 
the first series.  

Suppose the next `letter' $a(18)=2$.
Since $d^{(1)}_2<d^{(1)}_1$, that is, $q^{(1)}_1>0$, we
have a departure from the second queue in the first series,
that is, we decrease $q^{(1)}_1$ by one and increase $d^{(1)}_2$
by one.  In turn, this leads to an increase in $t^{(1)}_2$,
that is, a service at the second queue in the second series,
and so, since $q^{(2)}_1>0$, we decrease $q^{(2)}_1$ by one,
and increase $d^{(2)}_2$ by one.  That's it, the resulting
tableau $\tau(18)$ is
$$\begin{Young}
1&1&1&1&1&2&2&3&3\cr
2&2&2&2&3&3&3\cr
3&3\cr
\end{Young}$$
and we recognise this procedure as column-insertion of the
number 2 into the tableau $\tau(17)$.  

Note that, now, $q^{(2)}_1=0$.  Suppose the next letter $a(19)$ is
also a 2.  We still have $q^{(1)}_1>0$, so we decrease $q^{(1)}_1$ by one and increase $d^{(1)}_2$
by one.  In turn, this leads to an increase in $t^{(1)}_2$,
that is, a service at the second queue in the second series;
but $q^{(2)}_1=0$, so this service is unused and there is no
departure (that is, no increase in $d^{(2)}_2$).  The unused
service leads to an increase in $t^{(3)}_1$, that is, a service
at the only queue in the third series, and we increase $d^{(3)}_1$
by one.  The resulting tableau $\tau(19)$ is
$$\begin{Young}
1&1&1&1&1&2&2&3&3&3\cr
2&2&2&2&2&3&3\cr
3&3\cr
\end{Young}$$
and, again, we recognise this procedure as column-insertion of the
number 2 into the tableau $\tau(18)$.  And so on.

\hfill $\Box$

We will give a completely worked example, starting from an
empty tableau, in the next section.

We conclude this section with some remarks on the immediate
implications of Theorem~\ref{tgkrs}.  Let $l(n)$ and $\alpha(n)$
respectively denote the shape and weight of $\tau(n)$. 
In view of Theorem~\ref{tgkrs}, Lemma~\ref{prop4} can
now be interpreted as stating that $l_1(n)$ is the
length of the longest non-decreasing subsequence 
in the reversed word $a_n\cdots a_1$.  This is a well-known
property of the Robinson-Schensted algorithm.  Thus, Lemma~\ref{prop4} can
be regarded as a corollary of Theorem~\ref{tgkrs}, or the
proof of Lemma~\ref{prop4} given in the appendix
can be regarded as a new proof of the longest increasing
subsequence property of the Robinson-Schensted algorithm.

More generally, we can compare the statement of Theorem~\ref{tgkrs}
with Greene's theorem, and this leads to some remarkable identities.
Greene's theorem (see, for example,~\cite{llt})
states that, if $m_i(n)$ denotes the maximum of
the sum of the lengths of $i$ disjoint, non-decreasing subsequences
in the reversed word $a_n\cdots a_1$, then, for $i\le k$,
\begin{equation}
m_i(n)=l_1(n)+l_2(n)+\cdots+l_i(n).
\end{equation}
It therefore follows from Theorem~\ref{tgkrs} that
\begin{equation}\label{ri}
m_i(n)=G^{(k)}_k(n)+G^{(k)}_{k-1}(n)+\cdots+G^{(k)}_{k-i+1}(n).
\end{equation}
It would be interesting to see a direct proof of this identity.
Similarly, one can compare Theorem~\ref{tgkrs} with the various
extensions of Greene's theorem given, for example, in~\cite{bk}.

The implications of Lemma~\ref{prop2} for the Robinson-Schensted algorithm
would appear to be less well-known.
Fix $k\ge 2$, and set $H(z)=F^{(k)}(z^*)^*$.
\begin{cor} \label{recover} 
The weight $\alpha(n)$ of $\tau(n)$ can be 
recovered from the sequence of shapes $$\{\sh\tau(l),\ l\ge n\}.$$
In fact, if we set, for $m\ge 0$,
$$u(m)=\sh\tau(n+m)-\sh\tau(n) ,$$
then $$\alpha(n)=\sh\tau(n)+H(u).$$
\end{cor}

Recall that the recording tableaux $\sigma(n)$ are nested; the limiting
standard tableau $\sigma(\infty)=\lim_{n\to\infty}\sigma(n)$ is thus a well-defined 
object.  It follows from Corollary~\ref{recover} that we can recover the 
ininite word $a_1a_2\ldots$ from $\sigma(\infty)$.

This is also true for the Robinson-Schensted algorithm applied with
row insertion.  To see this, recall that the recording tableaux maintained
when one applies the row insertion algorithm to the infinite word $a_1a_2\ldots$
is the same as that maintained when one applies the column insertion algorithm
to the word $a_1^\dagger a_2^\dagger \ldots$, where $a_n^\dagger =k-a_n+1$ 
(or see, for example,~\cite[A.2, Exercise 5]{fulton}).  
We can thus extend the Robinson-Schensted correspondence
to a bijective mapping between infinite words and infinite standard tableaux.

\section{A worked example}

Suppose $k=3$ and $n=7$, and we apply
the Robinson-Schensted algorithm with column insertion to the word
$a_1\cdots a_7=3112322$.  
We obtain the following sequence of semistandard
tableaux:
$$
\begin{Young}
3 \cr
\end{Young}
\hspace{.5cm}
\begin{Young}
1& 3\cr
\end{Young}
\hspace{.5cm}
\begin{Young}
1&1&3\cr
\end{Young}
\hspace{.5cm}
\begin{Young}
1&1&3\cr
2\cr
\end{Young}
\hspace{.5cm}
\begin{Young}
1&1&3\cr
2\cr
3\cr
\end{Young}
\hspace{.5cm}
\begin{Young}
1&1&3\cr
2&2\cr
3\cr
\end{Young}
\hspace{.5cm}
\begin{Young}
1&1&2&3\cr
2&2\cr
3\cr
\end{Young}
$$

The evolution of the corresponding queueing network 
is as follows.  For each $n$, set
$$Q(n)=\left[
\begin{array}{cc}
q^{(1)}_1(n) & q^{(1)}_2(n) \\
& q^{(2)}_1(n) 
\end{array}\right]$$
and
$$D(n)=\left[
\begin{array}{ccc}
d^{(1)}_1(n) & d^{(1)}_2(n) & d^{(1)}_3(n) \\
& d^{(2)}_1(n) & d^{(2)}_2(n) \\
& & d^{(3)}_1(n) 
\end{array}\right] .$$
Initially,
$$Q(0)=\left[
\begin{array}{cc}
0 & 0\\
&0
\end{array}\right]
\ \ \ \ \mbox{and}\ \ \ \ 
D(0)=\left[\begin{array}{ccc}
0 & 0 & 0 \\
&0  & 0 \\
& &0 
\end{array}\right] .$$

At time 1, $a_1=3$, so
there is a service at the third queue in the first
series.  Since $q^{(1)}_2(0)=0$, this queue is empty, 
and the service is unused, leading to an increase
in $t^{(1)}_2$ and hence a service at the second 
queue in the second series.  This is also unused,
so we have a service at the first (and only) queue
in the third series, and there is a departure there:
$d^{(3)}_1(1)=1$.  Thus,
$$Q(1)=\left[
\begin{array}{cc}
0 & 0\\
&0
\end{array}\right]
\ \ \ \ \mbox{and}\ \ \ \
D(1)=\left[\begin{array}{ccc}
0 & 0 & 0 \\
&0  & 0 \\
& &1 
\end{array}\right] .$$
The corresponding tableau $\tau(1)$ is $$\begin{Young}3\cr\end{Young}$$

At time 2, $a_2=1$, so there is a service at the first queue in the first
series and a customer departs to join the second queue: thus,
$d^{(1)}_1(2)=1$ and $q^{(1)}_1(2)=1$.  In the second series,
this leads to an increase in $t^{(1)}_1$ and hence a departure
from the first queue to the second queue: thus, $d^{(2)}_1(2)=1$
and $q^{(2)}_1(2)=1$.  In turn, this leads to an increase in
$t^{(2)}_1$ and hence a service at the only queue in the third
series, which yields a (second) departure from that queue and
we have $d^{(3)}_1(2)=2$.
Thus,$$Q(2)=\left[
\begin{array}{cc}
1 & 0\\
&1
\end{array}\right]
\ \ \ \ \mbox{and}\ \ \ \
D(2)=\left[\begin{array}{ccc}
1 & 0 & 0 \\
&1  & 0 \\
& &2 
\end{array}\right] .$$
The corresponding tableau $\tau(2)$ is 
$$\begin{Young} 1&3\cr \end{Young}$$

Similarly, at time 3, $a_3=1$, and we get
$$Q(3)=\left[
\begin{array}{cc}
2 & 0\\
&2
\end{array}\right]
\ \ \ \ \mbox{and}\ \ \ \
D(3)=\left[\begin{array}{ccc}
2 & 0 & 0 \\
&2  & 0 \\
& &3
\end{array}\right] .$$
The corresponding tableau $\tau(3)$ is 
$$\begin{Young} 1&1&3\cr\end{Young}$$

At time 4, $a_4=2$, and there is a departure from the second
queue in the first series.  This provides a service at the 
second queue in the second series, which is non-empty, so we
have a departure from that queue as well.  Thus,
$$Q(4)=\left[
\begin{array}{cc}
1 & 1\\
&1
\end{array}\right]
\ \ \ \ \mbox{and}\ \ \ \
D(4)=\left[\begin{array}{ccc}
2 & 1 & 0 \\
&2  & 1 \\
& &3
\end{array}\right] .$$
The corresponding tableau $\tau(4)$
is $$\begin{Young} 1&1&3\cr 2\cr \end{Young}$$

At time 5, there is a service at the third queue in the first series,
and a customer departs.  That's all.
$$Q(5)=\left[
\begin{array}{cc}
1 & 0\\
&1
\end{array}\right]
\ \ \ \ \mbox{and}\ \ \ \
D(5)=\left[\begin{array}{ccc}
2 & 1 & 1 \\
&2  & 1 \\
& &3
\end{array}\right] .$$
The corresponding tableau $\tau(5)$
is $$\begin{Young} 1&1&3\cr 2\cr 3\cr \end{Young}$$

At time 6, there is a service at the second queue in the first series,
and a customer departs; this yields a service at the second queue in
the second series and a customer departs from there also.
$$Q(6)=\left[
\begin{array}{cc}
0 & 1\\
&0
\end{array}\right]
\ \ \ \ \mbox{and}\ \ \ \
D(6)=\left[\begin{array}{ccc}
2 & 2 & 1 \\
&2  & 2 \\
& &3
\end{array}\right] .$$
The corresponding tableau $\tau(6)$
is $$\begin{Young}  1&1&3\cr 2&2\cr 3\cr  \end{Young}$$

At time 7, there is a service at the second queue in the first series,
but this queue is empty and it is not used; this yields a service at
the {\em first} queue in the second series, leading to a departure
from that queue and consequently a service at (and departure from)
the first (and only) queue in the third series.  
$$Q(7)=\left[
\begin{array}{cc}
0 & 1\\
&1
\end{array}\right]
\ \ \ \ \mbox{and}\ \ \ \
D(7)=\left[\begin{array}{ccc}
2 & 2 & 1 \\
&3  & 2 \\
& &4
\end{array}\right] .$$
The final tableau $\tau(7)$ is 
$$\begin{Young} 1&1&2&3\cr 2&2\cr 3\cr \end{Young}$$

\section{Random walk in a Weyl chamber}

In this section we record some properties of the conditioned walk
of Theorem~\ref{oy2},
and extend its definition beyond the case $p_1<\cdots <p_k$.

The random walk $X$ is a Markov chain on $\Z_+^k$ with $X(0)=o$ 
and transition matrix
$$P(x,y) = p^{y-x} 1_{\{ y-x\in b\}} .$$
Denote by $P_x$ the law of the walk started at $x\in W\cap\Z_+^k$.

We will refer to the random walk with $p_1=\cdots=p_k$ 
as the {\em homogeneous walk}.

Denote by $s_l$ the Schur polynomial associated with
the integer partition $l_1\ge l_2\ge\cdots l_k\ge 0$.

\begin{lemma}
For any $r\in\P_k$, the function $h_r:\Z_+^k\to\R$, defined by
$$h_r(x)=p^{-x}s_{x^*}(r)1_{x\in W},$$ 
is harmonic for $P$.  Note that $h_r$ is strictly positive
on $W\cap\Z_+^k$.
\end{lemma}
{\em Proof.}  This follows immediately from the identity 
$$\sum_i s_{l+e_i}(r) = s_l(r),$$
which in turn can be seen as a special case of the Weyl character formula,
or verified directly using the formula
$$s_l(p)=\det\left(p_i^{l_j+k-j}\right) /\det\left(p_i^{k-j}\right) .$$ 
\hfill $\Box$

\begin{lemma}\label{ps}  Suppose $0<p_1<\cdots <p_k$.
Then, for any $x\in W\cap\Z_+^k$, 
$$P_x(X(n)\in W,\mbox{ for all }n\ge 0) = Cp^{-x}s_{x^*}(p),$$
where $C$ is a constant independent of $x$.  In particular,
the transition matrix associated with the conditioned walk
of Theorem~\ref{oy2} is given by
\begin{equation}\label{phat}
\hat P(x,y)=\frac{p^{-y}s_{y^*}(p)}{p^{-x}s_{x^*}(p)}P(x,y)
=\frac{s_{y^*}(p)}{s_{x^*}(p)}1_{\{ y-x\in b\}} .
\end{equation}
\end{lemma}
{\em Proof.}
When $r=p$, the Doob transform of the random walk $X$ via
the harmonic function $h_r$ has transition matrix $\hat P$.
Note that this can also be regarded as the Doob transform
of the homogeneous walk via the function $x\mapsto s_{x^*}(kp)$.
It follows from the asymptotic analysis of the Green function
associated with the (Poissonized) homogeneous walk presented
in~\cite{kor} that, if $\kappa(x,y)$ is the Martin kernel
associated with the homogeneous walk, then
$$\kappa(x,y)\to \mbox{\rm constant}\times s_{x^*}(kp)$$
whenever $y$ tends to infinity in $W$ in the direction $p$.
Thus, by standard Doob-Hunt theory (see, for example,~\cite{doob,w2}),
any realisation of the corresponding Doob transform, starting from 
the origin $o$, almost surely goes to infinity in the direction $p$. 
Moreover, any Doob transform on $W$ which almost surely goes to
infinity in the direction $p$ is necessarily the same Doob transform.
It therefore suffices to show that the `properly' conditioned walk of 
Theorem~\ref{oy2}, which is the Doob transform of $X$ via the harmonic 
function
$$g(x)=P_x(X(n)\in W,\mbox{ for all }n\ge 0) ,$$
almost surely goes to infinity in the direction $p$.  
But this follows immediately from the estimate,
denoting the law of the conditioned walk by $\hat P$,
$$\hat P( |X(n)-pn|>\epsilon n ) \le 
P( |X(n)-pn|>\epsilon n ) /g(x) \le Ke^{-c(\epsilon)n} /g(x),$$
where $c(\epsilon)>0$, and a standard Borel-Cantelli argument.
\hfill $\Box$

Note that the transition matrix $\hat P$ is well-defined by (\ref{phat})
for any $p\in\P_k$ and, by the symmetry of the Schur polynomials, is
symmetric in the $p_i$.  

The proof of Lemma~\ref{ps} given above is presented in more detail 
in~\cite{jphys}, where an explicit formula for the constant $C$ is
also given.

\section{The Robinson-Schensted algorithm with random words}\label{rsrw}

Having made the connection between the path-transformation
$G^{(k)}$ and the Robinson-Schensted algorithm, we will now 
give a direct proof of Theorem~\ref{oy2}, purely in the latter 
context. In fact, we will present a more general result,
which does not require the condition $p_1<\cdots <p_k$.

Let $\xi_1,\xi_2,\ldots$ be a sequence of independent random
variables with common distribution $p\in\P_k$.
Let $(S(n),T(n))$ be the pair of semistandard and standard tableaux
associated, by the Robinson-Schensted correspondence 
(with column-insertion), with the random word $\xi_1\xi_2\cdots\xi_n$. 
Here, $T(n)$ is the recording tableau.
Denote the shape of $S(n)$ by $\lambda(n)$; 
the weight (or type) of $S(n)$ is $X(n)$, where
$$X_i(n)=|\{1\le m\le n:\ \xi_m=i\}| .$$
Note that $X$ is the random walk discussed throughout this paper,
with transition matrix $$P(x,y) = p^{y-x} 1_{\{ y-x\in b\}} .$$

The joint law of $(S(n),T(n))$ is given, 
for $\sh\sigma=\sh\tau\vdash n$, by
\begin{equation}\label{jl}
P(S(n)=\sigma,T(n)=\tau)=p^\sigma ,
\end{equation}
and, for $x\in C=\{ x\in\Z_+^k:\ x_1\ge \cdots \ge x_k\}$,  
$$P(\lambda(n)=x)=s_x(p) f_x.$$
Here $p^\sigma$ is shorthand for $p^a$, where $a$ is the
weight of the tableau $\sigma$, and $f_x$ is the number
of standard tableaux with shape $x$.  The formula~(\ref{jl})
follows immediately from the fact that the Robinson-Schensted 
correspendence with column-insertion,
as in the case with row-insertion, is bijective.

Consider the Doob transform of $P$ on $C$, defined by its transition matrix
\begin{equation}
Q(x,y) = \frac{p^{-y}s_y(p)}{p^{-x}s_x(p)} P(x,y)
= \frac{s_y(p)}{s_x(p)}1_{\{ y-x\in b\}} .
\end{equation}

\begin{theorem}\label{dd} $\lambda$ is a Markov chain on $C$
with transition matrix $Q$.
\end{theorem}
{\em Proof.}  For $x,y\in C$, we will write $x\nearrow y$ if
$y-x\in b$.  Recall that a standard tableau $\tau$ with entries
$\{1,2,\ldots,n\}$ can be identified with a sequence of integer 
partitions
$$l(1)\nearrow l(2)\nearrow\cdots\nearrow l(n)$$
where $l(m)$ is the shape of the subtableau of $\tau$ consisting
only of the entries $\{1,2,\ldots,m\}$.  Since $T(n)$ is a recording
tableau, it is identified in this way with the sequence
$$\l(1)\nearrow\l(2)\nearrow\cdots\nearrow\l(n).$$
Thus, summing (\ref{jl}) over semistandard tableaux $\sigma$
with a given shape $l(n)\vdash n$, we obtain
\begin{equation}
P(\l(1)=l(1),\ldots,\l(n)=l(n)) = \sum_{\sh\sigma=l(n)} p^\sigma 
= s_{l(n)}(p) ,
\end{equation}
and so, for $x\nearrow y \vdash n+1$,
\begin{eqnarray*}
\lefteqn{ P(\l(n+1)=y | \l(1)=l(1),\ldots,\l(n-1)=l(n-1),\l(n)=x) }\\
\ \\
&=& \frac{ P(\l(1)=l(1),\ldots,\l(n-1)=l(n-1),\l(n)=x,\l(n+1)=y) }
{P(\l(1)=l(1),\ldots,\l(n-1)=l(n-1),\l(n)=x) }\\
&=& \frac{s_y(p)}{s_x(p)} ,
\end{eqnarray*}
as required.
\hfill $\Box$

Recalling the connection between $G^{(k)}$ and the Robinson-Schensted 
algorithm described in the previous section, and comparing $Q$ with the
transition matrix $\hat P$ defined by (\ref{phat}), we deduce 
the following generalisation of Theorem~\ref{oy2}.
\begin{cor}\label{ext}
$\hat X=G^{(k)}(X)$ is a Markov chain on $W\cap\Z_+^k$
with transition matrix $\hat P$.
\end{cor}

We will now record two lemmas which will yield an explicit description of
the joint law of $X$ and $\hat X$, and an intertwining relationship between 
their respective transition matrices.

Denote by $\kappa_{xy}$ the number of tableaux of shape $x$ and weight $y$. 
(These are the Kostka numbers.)

\begin{lemma}\label{wcfK}
\begin{equation}\label{wcf}
P(X(n)=y|\ \lambda(m),\ m\le n) = K(\lambda(n),y),
\end{equation}
where
\begin{equation}\label{K}
K(x,y) = \frac{p^y}{s_x(p)} \kappa_{xy} .
\end{equation}
\end{lemma}
{\em Proof.}
First note that the $\sigma$-algebra generated by $\{\lambda(m),\ m\le n\}$
is precisely the same as the $\sigma$-algebra generated by $T(n)$.
Thus, the conditional law of $X(n)$, given $\{\lambda(m),\ m\le n\}$,
is the same as the conditional law of $X(n)$, given $T(n)$.
But this only depends on the shape, $\lambda(n)$, of $T(n)$, and
is given by
$$K(x,y):=P(X(n)=y|\ \lambda(n)=x)=\frac{p^{y}}
{s_x(p)}\kappa_{xy},$$
as required.
\hfill $\Box$

We will now show that $P$ and $Q$ are intertwined
via the Markov kernel $K$, that is:
\begin{lemma}\label{int}
$QK=KP$. \end{lemma}
{\em Proof.}
\begin{eqnarray*}
(QK)(x,z) &=&
\sum_{y\in C} Q(x,y)K(y,z) \\
&=& \sum_{y\in C} P(x,y) \frac{p^{-y}s_y(p)}{p^{-x}s_x(p)}
\frac{p^z}{s_y(p)} \kappa_{yz} \\
&=& \sum_{i} p_i p^x p^{-x-e_i} \frac{p^z}{s_x(p)} \kappa_{x+e_i,z}\\
&=& \frac{p^z}{s_x(p)} \sum_i \kappa_{x+e_i,z} .
\end{eqnarray*}
On the other hand,
\begin{eqnarray*}
(KP)(x,z) &=& \sum_y K(x,y)P(y,z) \\
&=& \sum_i K(x,z-e_i) p_i \\
&=& \sum_i p_i \frac{p^{z-e_i}}
{s_x(p)} \kappa_{x,z-e_i} \\
&=& \frac{p^z}{s_x(p)} \sum_i \kappa_{x,z-e_i} .
\end{eqnarray*}
The statement of the lemma now follows from the identity
$$\sum_i \kappa_{x+e_i,z} = \sum_i \kappa_{x,z-e_i} ;$$
to see that this holds, observe that, for any $q\in\R_+^k$,
$$\sum_z q^z \sum_i \kappa_{x+e_i,z} = |q| s_x(q)
=\sum_z q^z \sum_i \kappa_{x,z-e_i}.$$
\hfill $\Box$

\begin{cor}\label{extB}
If we set $J(x,y)=K(x^*,y)$, then
\begin{equation}
P(X(n)=y|\ \hat X(m),\ m\le n,\ \hat X(n)=x)=J(x,y)
\end{equation}
and $\hat PJ=JP$.
\end{cor}

For a discussion on the role of intertwining in the context of Pitman's
$2M-X$ theorem (the case $k=2$), see~\cite{rp}.

It is instructive to note that Theorem~\ref{dd} also follows from 
Lemmas~\ref{wcfK} and~\ref{int}, provided we can show that there is a class of 
functions of the form $K\varphi$, $\varphi: \Z_+^k\to\R$,
which separate probability distributions on $C$.  
Indeed, by Lemmas~\ref{wcfK} and~\ref{int},
\begin{eqnarray*}
\lefteqn{ E\left[ (K\varphi)(\l(n+1)) |\ \l(m), m\le n \right] } \\
&=& E\left[ E\left[ \varphi(X(n+1))|\l(n+1) \right] | \l(m), m\le n \right] \\
&=& E\left[ E\left[ \varphi(X(n+1))|\l(m), m\le n+1\right] | \l(m), m\le n \right] \\
&=& E\left[ \varphi(X(n+1))|\l(m), m\le n\right] \\
&=& \sum_y K(\l(n),y) E\left[ \varphi(X(n+1))| X(n)=y \right] \\
&=& \sum_y K(\l(n),y) \sum_z P(y,z)\varphi(z) \\
&=& [(KP)\varphi](\l(n)) \\
&=& [(QK)\varphi](\l(n)) \\
&=& [Q(K\varphi)](\l(n)) 
\end{eqnarray*}
which would imply that $\l$ is a Markov chain with transition matrix $Q$
if the functions $K\varphi$ were determining.  To find such a class
of functions, we recall that the matrix $$\{\kappa_{xy}, (x,y)\in C^2\}$$
is invertible (see, for example,~\cite{macdonald}).  Thus, if we set,
for $q\in\R_+^k$,
\begin{equation}
\varphi_q(y) = p^{-y} \sum_{z\in C} \kappa^{(-1)}_{yz} q^z s_z(p) 1_{\{y\in C\}} ,
\end{equation}
we have $(K\varphi_q)(x)=q^x$, and these functions are clearly determining.

By exactly the same arguments as those given in the proof
of Theorem~\ref{dd}, if $\mu(n)$ denotes the shape
of the tableau obtained by applying the Robinson-Schensted algorithm with
{\em row} insertion to the random word $\xi_1\cdots\xi_n$,
we obtain:
\begin{theorem}  $\mu$ is a Markov chain on $C$ with transition matrix $Q$.
\end{theorem}

\section{Poissonized version}

We will now define a continuous version of $G^{(k)}$,
and state Poissonized versions of Corollaries~\ref{ext}
and~\ref{extB}.
This is an interesting setting in its own right, as the conditioned
walk in this case is closely related
to the Charlier ensemble and process (see, for example,
\cite{j2,kor}), but more importantly it provides a convenient
framework in which to apply Donsker's theorem and obtain the
Brownian analogue of Corollary~\ref{ext}, as was presented in~\cite{oy2} 
in the case $p=(1/k,\ldots,1/k)$.  Moreover, given the connection
we have now made with the Robinson-Schensted algorithm, this continuous path
transformation can also be regarded as a continuous analogue of
the Robinson-Schensted algorithm (see section 10 for further remarks in this
direction).

Let $D_0(\R_+)$ denote the space of cadlag paths $f:\R_+\to\R$
with $f(0)=0$.  We will extend the definition of the operations
$\i-c$ and $\s-c$ to a continuous context.  For $f,g\in D_0(\R_+)$,
define $f\i-c g\in D_0(\R_+)$ and $f\s-c g\in D_0(\R_+)$ by
\begin{equation}
(f\i-c g)(t)=\inf_{0\le s\le t}[f(s)+g(t)-g(s)] ,
\end{equation} 
and
\begin{equation}
(f\s-c g)(t)=\sup_{0\le s\le t}[f(s)+g(t)-g(s)] .
\end{equation}
As in the discrete case, these operations are not associative:
unless otherwise deleniated by parentheses, the default order of operations is
from left to right; for example, when we write $f\i-c g\i-c h$, we mean $(f\i-c g)\i-c h$.

Define a sequence of mappings $\Gamma^{(k)}:D_0(\R_+)^k\to D_0(\R_+)^k$ by
\begin{equation}
\Gamma^{(2)}(f,g)=(f\i-c g, g\s-c f),
\end{equation} 
and, for $k>2$,
\begin{align}
\Gamma^{(k)}(f_1, &  \ldots,f_k)=(f_1\i-c f_2\i-c\cdots\i-c f_k,\nonumber\\
& \Gamma^{(k-1)}(f_2\s-c f_1,f_3\s-c (f_1\i-c f_2),\ldots,f_k\s-c (f_1\i-c
\cdots\i-c f_{k-1}))) .
\end{align}

Let $N=(N_1,\ldots,N_k)$ be a continuous-time random walk with
generator $Gf(x)=\sum_i\mu_i [f(x+e_i)-f(x)].$
Denote by $\R_x$ the law of $N$ started from $x$, and by $(R_t)$
the corresponding semigroup. 
For convenience, we will also denote the corresponding
transition kernel by $R_t(x,y)$.  
Set $p_i=\mu_i/|\mu|$, and
denote by $\S_x$ the law of the $h_p$-transform on $W$, started at $x$,
and the corresponding semigroup by $S_t$.  

Note that the embedded discrete-time random
walk in $N$ has the same law as $X$.  That is,
if $\tau_n=\inf\{t\ge 0:\ |N(t)|=n\}$ and $Y(n)=N(\tau_n)$,
then $Y$ is a random walk on $\Z_+^k$ with transition matrix $P$.

\begin{theorem}\label{Poisson}
The law of $M = \Gamma^{(k)}(N)$ under $\R_o$ is the same
as the law of $N$ under $\S_o$.
\end{theorem}
Moreover, if $J$ is the Markov kernel defined in the previous section,
then:
\begin{theorem}\label{p1}
\begin{equation}
\R_o(N(t)=y|\ M(s),\ s\le t) = J(M(t),y) 
\end{equation}
and for $t\ge 0$, $S_tJ=JR_t$. 
\end{theorem}

\section{Brownian motion in a Weyl chamber and random matrices}

Let $X$ be a standard Brownian motion in $\R^d$,
and let $\Prob_x$ denote the law of $X$ started at $x$.
Denote the corresponding semigroup and transition kernel
by $(P_t)$, and the natural filtration of $X$ by $(\F_t)$.

Recall that $$W=\{ x\in\R^d:\ x_1\le \cdots\le x_d\}$$ and denote by
$\tilde P_t$ the semigroup of the process killed at
the first exit time
$$T=\inf\{ t\ge 0:\ X(t)\notin W\} .$$
Define $\Q_x$, for $x\in W^\circ$, by
$$\left. \Q_x\right|_{\F_t} = \frac{h(X(t\wedge T))}
{h(x)}\cdot \left. \Prob_x\right|_{\F_t} ,$$
where $h$ is the Vandermonde function $h(x)=\prod_{i<j}(x_j-x_i)$.
Denote the corresponding semigroup by $Q_t$.

The measure $$\Q_o=\lim_{W^\circ\ni x\to 0}\Q_x$$ 
is well-defined, and can be interpreted as
the law of the eigenvalue-process associated with Hermitian
Brownian motion~\cite{dyson,grabiner}.  
The law of $X(1)$ under $\Q_o$ is the
familiar Gaussian Unitary Ensemble (GUE) of random matrix theory.

In~\cite{oy2} it was shown, by applying Donsker's theorem
in the context of Theorem~\ref{Poisson} with $\mu=(1,\ldots,1)$,
that:
\begin{theorem}\label{oy2bm}
The law of $\Gamma^{(d)}(X)$ under $\Prob_o$ is the same
as the law of $X$ under $\Q_o$.  
\end{theorem}
In particular\footnote{See remark 2(ii) in section 10 below}, 
$(X_d\s-c\cdots\s-c X_1)(1)$ has the same law as the 
largest eigenvalue of a $d\times d$ GUE random matrix; this had been 
observed earlier by Baryshnikov~\cite{bar} and by Gravner, Tracy and 
Widom~\cite{gtw}.  A similar representation was obtained in~\cite{bj}.

Here we record some additional properties of the process $R=\Gamma^{(d)}(X)$, 
and its relationship with $X$, which are inherited, in the same application 
of Donsker's theorem, from Theorem~\ref{p1}.
\begin{theorem} \label{b1}
\begin{equation}
\Prob_o ( X(t)\in dx |\ R(s),\ s\le t; R(t)=r)
= L(r,dx),
\end{equation}
where $L$ is characterised by
\begin{equation}
\int_{\R^k} e^{\lambda\cdot y} L(x,dy) = \frac{\det(e^{\lambda_ix_j})}
{ h(\lambda) h(x)} =: c_\lambda(x) .
\end{equation}
Also, for $t\ge 0$, $Q_t L = L P_t$.
\end{theorem}
The intertwining $Q_t L = L P_t$ can also be seen as a direct
consequence of the Harish-Chandra/Itzykson-Zuber formula~\cite{hc,iz}
for the Laplace transform of the conditional law of the diagonal of a GUE 
random matrix given its eigenvalues, using the fact that the
diagonal of a Hermitian Brownian motion evolves according to 
the semigroup $P_t$ and the eigenvalues evolve according to the 
semigroup $Q_t$.
It is also easily verified using the Karlin-MacGregor formula 
$$\tilde P_t(x,y)= \det( P_t(x_i,y_j) ).$$

We will now present analogous results for Brownian motion with drift.
Fix $\mu\in\R^d$, and denote by $\Prob_x^{(\mu)}$ the law of Brownian
motion in $\R^d$ with drift $\mu$.  Denote by $(P_t^{(\mu)})$ the 
corresponding semigroup and by $(\tilde P_t^{(\mu)})$ the semigroup
of the process killed at the first exit time $T$ of the Weyl chamber
$W$.  Define $h_\mu$ by
\begin{equation}
h_\mu(x)= e^{-\mu\cdot x} \det\left( e^{\mu_i x_j}\right) .
\end{equation}
It is easy to check directly that $h_\mu$ a positive harmonic
function for $\tilde P_t^{(\mu)}$.  Define
$$\left. \Q^{(\mu)}_x\right|_{\F_t} = \frac{h_\mu(X(t\wedge T))}
{h_\mu(x)}\cdot \left. \Prob^{(\mu)}_x\right|_{\F_t} .$$
Denote the corresponding semigroup by $Q^{(\mu)}_t$.  Recalling the
absolute continuity relationship
\begin{equation}\label{acr}
\tilde  P^{(\mu)}_t (x,y) = e^{\mu\cdot (y-x) -|\mu|^2t/2} \tilde P_t(x,y),
\end{equation}
we can write
\begin{equation}
Q^{(\mu)}_t(x,y) = \frac{h_\mu(y)}{h_\mu(x)} \tilde P^{(\mu)}_t (x,y) 
= \frac{ \det\left( e^{\mu_i y_j}\right) }{ \det\left( e^{\mu_i x_j}\right) }
e^{-|\mu|^2t/2} \tilde P_t(x,y) ,
\end{equation}
and we note that this is symmetric in the $\mu_i$.

It is easy to verify that measure 
$$\Q^{(\mu)}_o=\lim_{W^\circ\ni x\to 0}\Q^{(\mu)}_x$$ 
is well-defined.

Applying Donskers theorem in the context of Theorem~\ref{Poisson},
as in~\cite{oy2}, we obtain:
\begin{theorem}
The law of $\Gamma^{(d)}(X)$ under $\Prob^{(\mu)}_o$ is the same
as the law of $X$ under $\Q^{(\mu)}_o$.  
\end{theorem}
As in the discrete case, we remark that the law $\Q^{(\mu)}_o$ 
is symmetric in the drifts $\mu_i$.

We also have, by the same application of Donsker's theorem,
the following analogue of Theorem~\ref{p1}.
\begin{theorem}
\begin{equation}
\Prob^{(\mu)}_o ( X(t)\in dx |\ R(s),\ s\le t; R(t)=r)
= L^{(\mu)}(r,dx),
\end{equation}
where 
\begin{equation}
L^{(\mu)}(x,dy) = c_\mu(x)^{-1} e^{\mu\cdot y} L(x,dy).
\end{equation}
Also, for $t\ge 0$, 
\begin{equation}\label{ird}
Q^{(\mu)}_t L^{(\mu)} = L^{(\mu)} P^{(\mu)}_t.
\end{equation}
\end{theorem}
The intertwining relationship (\ref{ird}) can also be verified directly 
using $Q_tL=L P_t$ and (\ref{acr}).

For related work on reflecting Brownian motions and non-colliding diffusions 
see~\cite{b0,bj,burdzy,bn,doumerc,hw,oy1,stw} and references therein.

\section{An application in queueing theory}

In this section, using the connection with the Robinson-Schensted correspondence
obtained in Section~\ref{gkrs}, we will write down a formula
for the `transient distribution' of a series of $M/M/1$ queues
in tandem.  There are many papers on this topic for the case
of a single queue, where the solution is given in terms of
modified Bessel functions; see~\cite{bm1} and references therein.
In~\cite{bm2}, the case of two queues was considered and a 
solution obtained, but the techniques used there do not seem
to extend easily to higher dimensions.  

Consider a series of $M/M/1$ queues in tandem, $k$ in number,
driven by Poisson processes $N_1, \ldots ,N_k$ with respective intensities
$\mu_1,\ldots,\mu_k$.  The first queue has infinitely many
customers, and the remaining queues are initially empty.
At every point of $N_i$, there is a service at the $i^{th}$
queue and, provided that queue is not empty, a customer
departs and joins the $(i+1)^{th}$ queue (or leaves the system
if $i=k$).

Denote by $D(t)=(D_1(t),\ldots,D_k(t))$ the respective numbers
of customers to depart from each queue up to time $t$.
Note that, since there are always infinitely many customers
in the first queue, $D_1$ is a Poisson process; we can
thus ignore the first queue, think of the second queue
as the first in a series of $k-1$ queues, and think of
$D_1$ as the arrival process at the first of these $k-1$ 
queues in series.  This is a more conventional set up
in queueing theory.  The state of the system is described
by the queue-lengths 
\begin{equation}
Q_1=D_1-D_2,\ \ldots\ ,\ Q_{k-1}=D_{k-1}-D_k .
\end{equation}
We will write down a formula for the law of $D(t)$
which, in turn, yields the law of $Q(t)=(Q_1(t),\ldots,Q_{k-1}(t))$.

Without loss of generality we can assume that $|\mu|=1$.
The de-Poissonized version of this problem is to consider
the usual random walk $X$, with $p=\mu$, and consider
the law of $\delta(n) = (D^{(k)}(X))(n)$.  
But we know this law, from sections~\ref{gkrs} and~\ref{rsrw}.
It is the law of $\beta(S(n))$, where $S(n)$ of the random semistandard tableau
obtained when one applies the Robinson-Schensted algorithm with column-insertion
to the random word $\xi_1\cdots\xi_k$ and $\beta_i(\tau)$
denotes the number of $i$'s in the $i^{th}$ row of a tableau $\tau$.

Thus, \begin{equation}\label{qf}
P(D(t)=d)=e^{-t}\sum_{n\ge 0} \frac{t^n}{n!} P(\delta(n) =d) ,
\end{equation}
where
\begin{equation}
P(\delta(n) =d) = \sum_{l\ge d, l\vdash n} \sum_{\sh\tau =l } p^\tau f_l 1_{\{\beta(\tau)=d\}} .
\end{equation}

This formula is complicated in general but simplifies in certain cases.

Consider the case $k=2$.  In this case, we have only one summand:
\begin{equation}
P(\delta(n)=d) = p_1^{d_1}p_2^{n-d_1} f_{(n-d_2,d_2)} .
\end{equation}
By the hook-length formula (see, for example,~\cite{fulton}), for $n\ge d_1+d_2$,
$$f_{(n-d_2,d_2)}= n! \frac{n-2d_2+1}{(n-d_2+1)! d_2!}.$$
Thus, recalling that $p=\mu$,
\begin{equation}
P(D(t)=d) =e^{-t}\sum_{n\ge d_1+d_2} t^n
\mu_1^{d_1}\mu_2^{n-d_1} \frac{n-2d_2+1}{(n-d_2+1)! d_2!} .
\end{equation}
It follows that
\begin{equation}
P(Q(t)=q)=(\mu_1/\mu_2)^q e^{-t} \sum_{m\ge q} (m+1) (\mu_2t)^m I_{m+1}
(2\sqrt{\mu_1\mu_2} t).
\end{equation}

Let $s_{l/d}$ denote the Schur polynomial associated with
the skew-tableau $l/d$ (see, for example,~\cite{fulton}).
We will use the following formula:
\begin{equation}\label{ss}
\sum_l s_{l/d}(x) \frac{t^{|l|}f_l}{|l|!} =  
e^{|x|t} \frac{t^{|d|}f_d}{|d|!} .
\end{equation}
This follows from the identity 
\begin{equation}
\sum_l s_{l/d}(x)s_l(y)=s_d(y)\prod_{i,j}(1-x_iy_j)^{-1}
\end{equation}
(this is a variant of Cauchy's identity; 
see, for example,~\cite[pp62-70]{macdonald}) and the fact that
\begin{equation}
\lim_{n\to\infty} s_l\left(\frac{t}{n} . 1^n\right) = \frac{t^{|l|}f_l}{|l|!} .
\end{equation}

In the case  $k=3$, if $p_2=p_3$, we have
\begin{eqnarray}\label{k3}
P(D_1(t)=d_1,D_2\ge d_2,D_3(t)=d_3) 
&=& e^{-t}p^d\sum_l s_{l/d}(p_2,p_3) \frac{t^{|l|}f_l}{|l|!}\nonumber \\
&=& e^{-p_1t} p^d  \frac{t^{|d|}f_d}{|d|!} .
\end{eqnarray}
It would be interesting to compare (\ref{k3}) with the explicit formulas 
obtained in~\cite{bm2} for this case.  

In the general case, we can simplify the formula (\ref{qf}) if $d$ is constant.
Suppose $d_i=m$ for all $i$.  Then, using~(\ref{ss}) and the hook-length formula,
\begin{eqnarray}
P(D(t)=d)&=&
e^{-t}p^d\sum_l s_{l/d}(p_2,p_3,\ldots,p_k) \frac{t^{|l|}f_l}{|l|!} \\
&=& e^{-p_1t} p^d \frac{t^{|d|}f_d}{|d|!} \\
&=& e^{-p_1t}\prod_ip_i^m t^{mk} \frac{G(k+l)}{G(k)G(l)} ,
\end{eqnarray}
where $G$ is Barne's function, defined by $G(k)=\prod_{i\le k} \Gamma(i)$.  

\section{Concluding remarks}

{\em 1. The Krawchouk process:}
A binomial version of Theorem~\ref{oy2} was presented in \cite{kor}.
This states that, if $X$ is a random walk in $Z_+^k$ with transition
matrix $$P(x,y)=Cp^{y-x}1_{y-x\in\{0,1\}^k}$$
($C$ is a normalising constant) then, 
assuming $p_1<\cdots <p_k$ and extending slightly the domain of $G^{(k)}$,
$G^{(k)}(X)$ has the same
law as that of $X$ conditioned to stay forever in $W$.
In this case, a similar connection can be made with the dual Robinson-Schensted-Knuth (RSK)
correspondence for zero-one matrices and analogues of all the main results 
of Section~\ref{rsrw} can be obtained similarly.  The Schur polynomials again play an important role.
See~\cite{jphys} for details.

{\em 2. Properties of $\Gamma^{(k)}$:}
The following continuous analogues of Lemmas~\ref{prop1}-\ref{prop4} 
can be readily verified.  Denote by $C_0(\R_+,\R^k)$ the set of continuous
functions $f:\R_+\to\R^k$ with $f(0)=o$.  For $f\in C_0(\R_+,\R^k)$, where $k\ge 2$,
\begin{itemize}
\item[(i)] $|\Gamma^{(k)}(f)|=|f| $
\item[(ii)] $\Gamma_k^{(k)}(f)=f_k\s-c\cdots\s-c f_1$
\item[(iii)] $f(t)=[\Gamma^{(k)}(f)](t)
+\Phi^{(k)}([\Gamma^{(k)}(f)](t,u),\ u\ge t),$
where $\Phi^{(k)}$ is defined on a suitable domain as the
continuous analogue of $F^{(k)}$.
\end{itemize}
In this continuous setting, the identity (ii) can be verified directly
using the `sup-integration by parts' formula
$$\sup_{0<s<t} \left\{ \sup_{0<r<s} u(r) + v(s) \right\}
\bigvee \sup_{0<s<t} \left\{ u(s)+ \sup_{0<r<s} v(r) \right\} 
\hspace{1in}$$
$$\hspace{3in}  =\sup_{0<s<t} u(s)+ \sup_{0<s<t} v(s) ,$$
for $u,v\in C_0(\R_+,\R^k)$.  This is a `max-plus' analogue of the
usual integration by parts formula and is easily verified by the
method of Laplace.

{\em 3. A continuous Robinson-Schensted algorithm:}
Given the connection we have made between the path-transformations
$G^{(k)}$ and the Robinson-Schensted algorithm, the mappings $\Gamma^{(k)}$ can be
used to define a continuous version of the Robinson-Schensted algorithm.  More precisely,
let $\C_{GC}$ denote the {\em Gelfand-Cetlin cone}, which consists of
triangular arrays of real numbers
\begin{equation}
(x^{(i)}_j,\ 1\le i\le k,\ 1\le j\le i)
\end{equation}
satisfying $x^{(i)}_j\ge x^{(i-1)}_j \ge x^{(i)}_{j+1}$, for all $i,j$.
Points in the Gelfand-Cetlin cone can be regarded as continuous
analogues of semistandard tableaux.  The continuous analogue of a
word is a continuous function $f:[0,1]\to \R^k$ with
$f(0)=o$: denote the set of these functions by $C_0([0,1],\R^k)$.
Define a map $\phi:C_0([0,1],\R^k)\to \C_{GC}$, as follows.
For convenience, let $\Gamma^{(1)}$ be the identity transformation.
If we set $x=\phi(f_1,\ldots,f_k)$ then, for each $1\le i\le k$,
\begin{equation}
x^{(i)}=([\Gamma_i^{(i)}(f_1,\ldots,f_i)](1),\ldots,
[\Gamma_1^{(i)}(f_1,\ldots,f_i)](1)).
\end{equation}
The continuous analogue of the corresponding `recording tableau'
is the path
\begin{equation}
\rho(f)= \{ [\Gamma^{(k)}(f)](t) ,\ 0\le t\le 1\} \in C_0([0,1],W) .
\end{equation}
By analogy with the discrete Robinson-Schensted algorithm, the function $f$ can
be uniquely recovered from the pair $\phi(f)$ and $\rho(f)$.
A more detailed discussion on the properties of this continuous Robinson-Schensted 
algorithm will be presented elsewhere.

{\em 4. GUE minors:}
Let $A$ be a $k\times k$ GUE random matrix and denote the
eigenvalues of the $i^{th}$ minor $(A_{lm},\ l,m\le i)$
by $\lambda_1^{(i)}\ge \cdots \ge \lambda_i^{(i)}$, for $i\le k$.
In the above context, Baryshnikov~\cite{bar} showed that,
if $(B(t),\ 0\le t\le 1)$ is a standard Brownian motion
in $\R^k$, then the random vector $$(\phi_1^{(1)}(B),\ldots,
\phi_1^{(k)}(B))$$ has the same law as $$(\lambda_1^{(1)},\ldots,
\lambda_1^{(k)}).$$
In~\cite{bar}, Donsker's theorem is applied in the context
of a random semistandard tableau with the same law as $T(n)$
of section 6 in the homogeneous case $p_1=\cdots=p_k$.
We can thus extend Baryshnikov's arguments using the
the representation for the Robinson-Schensted algorithm given in this
paper and the continuity of the mappings $\Gamma^{(i)}$,
$i\le k$, to see that, in fact,
$\phi(B)$ has the same law as $$(\lambda^{(1)},\ldots,
\lambda^{(k)}).$$
However, it is easy to 
see\footnote{Bougerol and Jeulin, private communication}, 
by considering the case $k=2$,
that these identities do not extend to the process level
(that is, with $\phi$ defined simultaneously on intervals
$[0,t]$ instead of just $[0,1]$ and `GUE' replaced by
`Hermitian Brownian motion').

{\em 5. Related topics:}
The intertwining of Lemma~\ref{int} is closely related to
the work of Biane on quantum random walks~\cite{b1,b2,b3}.
The Robinson-Schensted correspondence is of fundamental importance to the
representation theory of $S_n$ and $GL(n)$.
Related topics in representation theory (which are certainly connected
to results presented in this paper) include Littelmann's path model 
for the finite-dimensional representations of $GL(n)$ 
(see, for example,~\cite{littelmann}), crystal bases and 
representations of quantum groups (see, for example,~\cite{djm,llt}).  

\section*{Appendix}

{\em Proof of Lemma \ref{prop1}:}  The first identity is trivial:
\begin{eqnarray*}
(y\s-c x)(n) &=& \max_{0\le m\le n} [ y(m)+x(n)-x(m) ] \\
&=& x(n)+y(n)+\max_{0\le m\le n} [ y(m)-y(n)-x(m) ]\\
&=& x(n)+y(n)-(x\i-c y)(n) .
\end{eqnarray*}
If we set $z=y-x$, and $s(n)=\max_{0\le m\le n} z(n)$, then
the second identity is equivalent to the well-known fact that
$$s(n)=\min_{l\le n} [2s(l)-z(l)] .$$
\hfill $\Box$

{\em Proof of Lemma \ref{prop2}:}
Fix $x\in\Lambda_k$, and write $G^{(k)}=G^{(k)}(x)$ unless otherwise
indicated.  Similarly for $D^{(k)}$ and $T^{(k)}$.  We will first
show that $|G^{(k)}|=|x|$. 
We will prove this by induction on $k$.  The case $k=2$ is given
by Lemma~\ref{prop1}.  Assume the induction hypothesis
for $k-1$.
We recall from the definitions that
\begin{equation}
G^{(k)}=\left( D^{(k)}_k, G^{(k-1)}\left( T^{(k)}\right) \right) .
\end{equation}
By the induction hypothesis,
\begin{equation}
|G^{(k)}|=D^{(k)}_k+|T^{(k)}| .
\end{equation}
Recall that, for $i\ge 2$,
\begin{equation}
D_i^{(k)}=D^{(k)}_{i-1}\i-c x_i ,
\end{equation}
and 
\begin{equation}
T^{(k)}_{i-1}=x_i\s-c D^{(k)}_{i-1}.
\end{equation}
Thus, by Lemma~\ref{prop1},
\begin{equation}\label{xi}
x_i = D_i^{(k)}-D^{(k)}_{i-1}+T^{(k)}_{i-1},
\end{equation}
for each $i\ge 2$; summing this over $i$ and
recalling that $D^{(k)}_1=x_1$ yields
\begin{equation}
|x| = D_k^{(k)} + |T^{(k)}|,
\end{equation}
so we are done.

We will now show that
\begin{equation}
x(n)=G^{(k)}(n)+F^{(k)}\left( G^{(k)}(l)-G^{(k)}(n),\ l\ge n\right) ,
\end{equation}
for some function $F^{(k)}$ to be defined.
Again we will prove this by induction on $k$, and note that
for $k=2$, this is given by Lemma~\ref{prop1}.

A recursive definition of $F^{(k)}$ will be implicit in the
induction argument.  Recall that
\begin{equation}
G^{(k)}=\left( D^{(k)}_k, G^{(k-1)}\left( T^{(k)}\right) \right) .
\end{equation}
Assuming the induction hypothesis for $k-1$, we have, for $i\ge 2$,
\begin{equation}
T_{i-1}^{(k)}(n) = G_i^{(k)}(n)+F^{(k-1)}\left(
(G^{(k)}_2,\ldots,G^{(k)}_k)(n,l) \right) .
\end{equation}
Thus, for $i\ge 2$, using (\ref{xi}) and the fact that 
$$D_i^{(k)}(n)-D^{(k)}_{i-1}(n)=\max_{l\ge n}
[ D_i^{(k)}(n,l)-T_{i-1}^{(k)}(n,l) ],$$ we have
\begin{equation}\label{yy}
x_i(n) = G_i^{(k)}(n) + J^{(k)}_i\left(
(D_i^{(k)},G^{(k)}_2,\ldots,G^{(k)}_k)(n,l) \right) ,
\end{equation}
where $J^{(k)}_i$ is defined on a suitable domain.
It is important to note here is that $J^{(k)}_i$ does not
depend on $n$.  

In this way, recalling that $D^{(k)}_k=G^{(k)}_1$, we obtain
\begin{equation}
x_k(n) = G_k^{(k)}(n)+F_k^{(k)}\left( G^{(k)}(l)-G^{(k)}(n),\ l\ge n\right) ,
\end{equation}
where the function $F_k^{(k)}$ is implicitly defined by this identity
(on a suitable domain) and does not depend on $n$.
Observe that we can also recover the sequence of future increments 
$x_k(l)-x_k(n)$ as a function, which does not depend on $n$, of the
sequence $\{G^{(k)}(l)-G^{(k)}(n),\ l\ge n\}$.  

We can now recover the values $x_{k-1}(n)$, $x_{k-2}(n)$, and so on,
as follows.  By equations (\ref{xi}) with $i=k$,
\begin{equation}
D^{(k)}_{k-1}(n)=G_1^{(k)}(n)-x_k(n)+T_{k-1}^{(k)}(n).
\end{equation}
It follows that the sequence $\{D^{(k)}_{k-1}(n,l),\ l\ge n\}$ 
is a function, which does not depend on $n$, of the
sequence $\{G^{(k)}(l)-G^{(k)}(n),\ l\ge n\}$.  Combining this 
with~(\ref{yy}), we see that 
\begin{equation}
x_{k-1}(n) = 
G_{k-1}^{(k)}(n)+F_{k-1}^{(k)}\left( G^{(k)}(l)-G^{(k)}(n),\ l\ge n\right) ,
\end{equation}
where $F_{k-1}^{(k)}$ is implicitly defined by this identity
(on a suitable domain) and does not depend on $n$.
Similarly, we can recover the sequence of future increments
$x_k(l)-x_k(n)$ as a function, which does not depend on $n$, of the
sequence $\{G^{(k)}(l)-G^{(k)}(n),\ l\ge n\}$.  And so on.
Finally, $x_1(n)$ is obtained using $|x|=|G^{(k)}|$.
\hfill $\Box$

{\em Proof of Lemma \ref{prop3}:}
We want to show that, for $(a,b,c)\in\Lambda_3$,
\begin{equation}\label{uf}
a\s-c (c\i-c b)\s-c (b\s-c c) = a\s-c b \s-c c,
\end{equation}
and for $(w,x,y)\in\Lambda_3$,
\begin{equation}\label{ouf}
w\i-c (y\s-c x)\i-c (x\i-c y) = w\i-c x\i-c y .
\end{equation}

First note that these identities are equivalent.
To see this, set $a(n)=n-w(n)$, $b(n)=n-x(n)$ and
$c(n)=n-y(n)$, then plug these into (\ref{uf}) to obtain
(\ref{ouf}).  We will therefore restrict our attention
to the identity~(\ref{ouf}).

Let $d=x\i-c y$, $t=y\s-c x$, $q=x-d$ and $u=y-d$.
Then~(\ref{ouf}) becomes
\begin{equation}
w\i-c (x+u)\i-c (y-u) = w\i-c x\i-c y .
\end{equation}
That is, the output of a series of queues in tandem
driven by $(w,x+u,y-u)$ is the same as that of the
series driven by $(w,x,y)$.  Set 
\begin{eqnarray*}
d_1 &=& w\i-c x\\
d_2 &=& w\i-c x\i-c y \\
\tilde d_1 &=& w\i-c (x+u)\\
\tilde d_2 &=& w\i-c (x+u)\i-c (y-u)
\end{eqnarray*}
and
\begin{eqnarray*}
q_1 &=& w-d_1\\
q_2 &=& d_1-d_2\\
\tilde q_1 &=& w-\tilde d_1\\
\tilde q_2 &=& \tilde d_1 - \tilde d_2 .
\end{eqnarray*}
We want to show that $d_2=\tilde d_2$.
From the above definitions, this is equivalent to
showing that $$q_1(n)+q_2(n)=\tilde q_1(n)+\tilde q_2(n)$$
for all $n\ge 0$.  We will prove this by induction on $n$.

The induction hypothesis $H$ is:

$q_1+q_2=\tilde q_1+\tilde q_2$, and {\em either}

(i) $\tilde q_2-q_2\ge 0\mbox{ and }q-q_2=0$, {\em or}

(ii) $\tilde q_2-q_2=0\mbox{ and }q-q_2\ge 0$.

When $n=0$, $q=q_1=q_2=\tilde q_1=\tilde q_2=0$, and the
induction hypothesis is trivially satisfied.
Assume the induction hypothesis holds at time $n-1$.
Note that $(w,x,y-u,u)\in\Lambda_4$, that is, only
one of these quantities, if any, can increase by one
at time $n$.  We will consider the following five cases,
which are exaustive and mutually exclusive, separately.
\begin{itemize}
\item[(a)] $(w,x,y-u,u)(n)=(w,x,y-u,u)(n-1)$
\item[(b)] $w(n)-w(n-1)=1$
\item[(c)] $x(n)-x(n-1)=1$
\item[(d)] $(y-u)(n)-(y-u)(n-1)=1$
\item[(e)] $u(n)-u(n-1)=1$
\end{itemize}

{\bf Case (a):} $(w,x,y-u,u)(n)=(w,x,y-u,u)(n-1)$.  In this
case, nothing changes, and so $H$ is preserved.

{\bf Case (b):} $w(n)-w(n-1)=1$.
In this case, $q_1(n)=q_1(n-1)+1$ and $\tilde q_1(n)=\tilde q_1(n-1)+1$,
the other quantities remain unchanged, and $H$ is preserved.

{\bf Case (c):} $x(n)-x(n-1)=1$.  Then $q(n)=q(n-1)+1$.

Suppose, $q_1(n-1)>\tilde q_1(n-1)>0$.  Then $q_1(n)=q_1(n-1)-1$
and $q_2(n)=q_2(n-1)+1$.  Thus, $q-q_2$ and $q_1+q_2$ do not change.
Also, $\tilde q_1(n)=\tilde q_1(n-1)-1$ and $\tilde q_2(n)=\tilde q_2(n-1)+1$. 
Thus, $\tilde q_2-q_2$ and $\tilde q_1+\tilde q_2$ do not change either, 
so we still have $q_1+q_2=\tilde q_1+\tilde q_2$, and $H$ is preserved.

Now suppose $q_1(n-1)>\tilde q_1(n-1)=0$. Note that this implies
$\tilde q_2(n-1)-q_2(n-1)>0$, so that we are initially in case (i)
of the induction hypothesis. In this case,
$q_1(n)=q_1(n-1)-1$ and $q_2(n)=q_2(n-1)+1$, but $\tilde q_1$
and $\tilde q_2$ do not change.  Thus, $q-q_2$, $q_1+q_2$ and
$\tilde q_1+\tilde q_2$ do not change.  The quantity $\tilde q_2-q_2$
decreases by one, but remains non-negative, so we remain in case (i)
and $H$ is preserved.  

Finally, if $q_1(n-1)=\tilde q_1(n-1)=0$, then we are initially
in case (ii) of $H$.  There is no change to $q_1,q_2,\tilde q_1$
or $\tilde q_2$, but $q$ inceases by one and so we remain in case
(ii) and $H$ is preserved.

{\bf Case (d):} $(y-u)(n)-(y-u)(n-1)=1$.
In this case, $q(n-1)>0$ and $q$ decreases by one.
The values of $q_1$ and $\tilde q_1$ do not change.

If we are in case (i) at time $n-1$, then $\tilde q_2(n-1)\ge q_2(n-1)>0$
and so both $q_2$ and $\tilde q_2$ also decrease by one; thus, we
remain in case (i) and $H$ is preserved.

If we are in case (ii) at time $n-1$, then $\tilde q_2(n-1)= q_2(n-1)$
and either $q_2$ and $\tilde q_2$ both decrease by one or both
remain unchanged.  Either way, we remain in case (ii) and $H$ is
preserved.

{\bf Case (e):} $u(n)-u(n-1)=1$.  Then $q(n-1)=q_2(n-1)=0$ and
we are initially in case (i) of $H$.  The values of $q,q_1$ and
$q_2$ will not change.  If $\tilde q_1>0$, then $\tilde q_1$
decreases by one and $\tilde q_2$ increases by one; otherwise,
$\tilde q_1$ and $\tilde q_2$ do not change.  Either way, we
remain in case (i) and $H$ is preserved.
\hfill $\Box$

{\em Proof of Lemma \ref{prop4}:}
We will prove this by induction on $k$.  It is certainly
true for $k=2$, from the definition of $G^{(2)}$.  
Recalling the definition of $G^{(k)}$,
\begin{equation}
G^{(k)}=\left( D^{(k)}_k, G^{(k-1)}\left( T^{(k)}\right) \right) .
\end{equation}
By the induction hypothesis, that the Theorem is true for $G^{(k-1)}$,
\begin{equation}
G^{(k)}_k = T^{(k)}_{k-1}\s-c T^{(k)}_{k-2}\s-c \cdots T^{(k)}_1 .
\end{equation}
We will now repeatedly apply Lemma \ref{prop3}.
\begin{eqnarray*}
T^{(k)}_{k-1}\s-c T^{(k)}_{k-2} &=& x_k\s-c D^{(k)}_{k-1} \s-c T^{(k)}_{k-2}\\
&=& x_k\s-c (D^{(k)}_{k-2}\i-c x_{k-1})\s-c (x_{k-1}\s-c D^{(k)}_{k-2})\\
&=& x_k\s-c x_{k-1}\s-c D^{(k)}_{k-2} .
\end{eqnarray*}
Similarly,
\begin{eqnarray*}
x_k\s-c x_{k-1}\s-c D^{(k)}_{k-2}\s-c T^{(k)}_{k-3} &=& 
x_k\s-c x_{k-1}\s-c (D^{(k)}_{k-3}\i-c x_{k-2})\s-c (x_{k-2}\s-c D^{(k)}_{k-3})\\
&=& x_k\s-c x_{k-1}\s-c x_{k-2}\s-c D^{(k)}_{k-3} ,
\end{eqnarray*}
and so on.
\hfill $\Box$

\end{document}